\documentclass[12pt,reqno,openbib,runningheads,a4paper]{amsart}%
\usepackage{eurosym}
\usepackage{amssymb}
\usepackage{amsfonts}
\usepackage{amsmath}
\usepackage{graphicx}
\usepackage[authoryear]{natbib}
\usepackage[dvipsnames]{xcolor}
\usepackage[legalpaper,bookmarks=true,colorlinks=true,linkcolor=blue,citecolor=blue]%
{hyperref}%
\providecommand{\U}[1]{\protect\rule{.1in}{.1in}}
\providecommand{\U}[1]{\protect\rule{.1in}{.1in}}
\newtheorem{theorem}{Theorem}[section]

\makeatletter
\renewcommand{\@biblabel}[1]{}
\makeatother
\begin{document}

\begin{center}
{\Large \textbf{A Lynden-Bell integral estimator for the tail index of
right-truncated data with a random threshold}}\medskip

{\large Nawel Haouas, Abdelhakim Necir\footnote{\noindent{\small Corresponding
author: \texttt{necirabdelhakim@yahoo.fr} \newline}\noindent
{\small \textit{E-mail addresses:}\newline\texttt{nawel.haouas@yahoo.com}
(N.~Haouas)\newline\texttt{djmeraghni@yahoo.com} (D.~Meraghni)\newline%
\texttt{brah.brahim@gmail.com} (B.~Brahimi)}}, Djamel Meraghni, Brahim
Brahimi}\medskip

{\small \textit{Laboratory of Applied Mathematics, Mohamed Khider University,
Biskra, Algeria}}\medskip\medskip

\end{center}

\noindent\textbf{Abstract}\medskip

\noindent By means of a Lynden-Bell integral with deterministic threshold,
Worms and Worms [A Lynden-Bell integral estimator for extremes of randomly
truncated data. Statist. Probab. Lett. 2016; 109: 106-117] recently introduced
an asymptotically normal estimator of the tail index for randomly
right-truncated Pareto-type data. In this context, we consider the random
threshold case to derive a Hill-type estimator and establish its consistency
and asymptotic normality.\ A simulation study is carried out to evaluate the
finite sample behavior of the proposed estimator.\medskip


\medskip\medskip

\section{\textbf{Introduction\label{sec1}}}

\noindent Let $\left(  \mathbf{X}_{i},\mathbf{Y}_{i}\right)  ,$ $1\leq i\leq
N$ be a sample of size $N\geq1$ from a couple $\left(  \mathbf{X}%
,\mathbf{Y}\right)  $ of independent random variables (rv's) defined over some
probability space $\left(  \Omega,\mathcal{A},\mathbf{P}\right)  ,$ with
continuous marginal distribution functions (df's) $\mathbf{F}$ and
$\mathbf{G}$ respectively.$\ $Suppose that $\mathbf{X}$ is truncated to the
right by $\mathbf{Y},$ in the sense that $\mathbf{X}_{i}$ is only observed
when $\mathbf{X}_{i}\leq\mathbf{Y}_{i}.$ We assume that both survival
functions $\overline{\mathbf{F}}:=1-\mathbf{F}$\textbf{\ }and $\overline
{\mathbf{G}}:=1-\mathbf{G}$ are regularly varying at infinity with respective
negative indices $-1/\gamma_{1}$ and $-1/\gamma_{2}.$ That is, for any $x>0,$%
\begin{equation}
\lim_{z\rightarrow\infty}\frac{\overline{\mathbf{F}}\left(  xz\right)
}{\overline{\mathbf{F}}\left(  z\right)  }=x^{-1/\gamma_{1}}\text{ and }%
\lim_{z\rightarrow\infty}\frac{\overline{\mathbf{G}}\left(  xz\right)
}{\overline{\mathbf{G}}\left(  z\right)  }=x^{-1/\gamma_{2}}. \label{RV-1}%
\end{equation}
It is well known that, in extreme value analysis, weak approximations are
achieved in the second-order framework (see, e.g., \citeauthor{deHF06},
\citeyear[page 48]{deHF06}). Thus, it seems quite natural to suppose that
$\mathbf{F}$\ and $\mathbf{G}$\ satisfy the second-order condition of regular
variation, which we express in terms of the tail quantile functions pertaining
to both df's. That is, we assume that for $x>0,$ we have%
\begin{equation}
\underset{t\rightarrow\infty}{\lim}\dfrac{\mathbb{U}_{\mathbf{F}}\left(
tx\right)  /\mathbb{U}_{\mathbf{F}}\left(  t\right)  -x^{\gamma_{1}}%
}{\mathbf{A}_{\mathbf{F}}\left(  t\right)  }=x^{\gamma_{1}}\dfrac{x^{\tau_{1}%
}-1}{\tau_{1}}, \label{second-order}%
\end{equation}
and%
\begin{equation}
\underset{t\rightarrow\infty}{\lim}\dfrac{\mathbb{U}_{\mathbf{G}}\left(
tx\right)  /\mathbb{U}_{\mathbf{G}}\left(  t\right)  -x^{\gamma_{2}}%
}{\mathbf{A}_{\mathbf{G}}\left(  t\right)  }=x^{\gamma_{2}}\dfrac{x^{\tau_{2}%
}-1}{\tau_{2}}, \label{second-orderG}%
\end{equation}
where $\tau_{1},\tau_{2}<0$\ are the second-order parameters and
$\mathbf{A}_{\mathbf{F}},$ $\mathbf{A}_{\mathbf{G}}$\ are functions tending to
zero and not changing signs near infinity with regularly varying absolute
values at infinity with indices $\tau_{1},$ $\tau_{2}$ respectively. For any
df $K,$ the function $\mathbb{U}_{K}\left(  t\right)  :=K^{\leftarrow}\left(
1-1/t\right)  ,$ $t>1,$ stands for the tail quantile function, with
$K^{\leftarrow}\left(  u\right)  :=\inf\left\{  v:K\left(  v\right)  \geq
u\right\}  ,$ $0<u<1,$ denoting the generalized inverse of $K.$ From Lemma 3
in \cite{HJ-2011} , the second-order conditions $(\ref{second-order})$ and
$(\ref{second-orderG})$ imply that there exist constants $d_{1},d_{2}>0,$ such
that%
\begin{equation}
\overline{\mathbf{F}}\left(  x\right)  =d_{1}y^{-1/\gamma_{1}}\ell_{1}\left(
x\right)  \text{ and }\overline{\mathbf{G}}\left(  x\right)  =d_{2}%
y^{-1/\gamma_{2}}\ell_{2}\left(  x\right)  ,\text{ }x>0, \label{HJ}%
\end{equation}
where $\lim_{x\rightarrow\infty}\ell_{i}\left(  x\right)  =1$ and $\left\vert
1-\ell_{i}\right\vert $ is regularly varying at infinity with tail index
$\tau_{i}\gamma,$ $i=1,2.$ This ccondition is fullfilled by many commonly used
models such as Burr, Fr\'{e}chet, Generalized Pareto, absolute Student,
log-gamma distributions, to name but a few. Also known as heavy-tailed,
Pareto-type or Pareto-like distributions, these models take a prominent role
in extreme value theory and have important practical applications as they are
used rather systematically in certain branches of non-life insurance, as well
as in finance, telecommunications, hydrology, etc...
\citep[see, e.g.,][]{Res06}. $\medskip$

\noindent Let us now denote $\left(  X_{i},Y_{i}\right)  ,$ $i=1,...,n$ to be
the observed data, as copies of a couple of rv's $\left(  X,Y\right)  ,$
corresponding to the truncated sample $\left(  \mathbf{X}_{i},\mathbf{Y}%
_{i}\right)  ,$ $i=1,...,N,$ where $n=n_{N}$ is a sequence of discrete rv's
which, in virtue of the weak law of large numbers, satisfies $n_{N}%
/N\overset{\mathbf{P}}{\rightarrow}p:=\mathbf{P}\left(  \mathbf{X}%
\leq\mathbf{Y}\right)  ,$ as $N\rightarrow\infty.$ We denote the joint df of
$X$ and $Y$ by $H\left(  x,y\right)  :=\mathbf{P}\left(  X\leq x,Y\leq
y\right)  =\mathbf{P}\left(  \mathbf{X}\leq\min\left(  x,\mathbf{Y}\right)
,\mathbf{Y}\leq y\mid\mathbf{X}\leq\mathbf{Y}\right)  ,$ which is equal to
$p^{-1}%
{\displaystyle\int_{0}^{y}}
\mathbf{F}\left(  \min\left(  x,z\right)  \right)  d\mathbf{G}\left(
z\right)  .$ The marginal distributions of the rv's $X$ and $Y,$ respectively
denoted by $F$ and $G,$ are given by $F\left(  x\right)  =p^{-1}%
{\displaystyle\int_{0}^{x}}
\overline{\mathbf{G}}\left(  z\right)  d\mathbf{F}\left(  z\right)  $ and
$G\left(  y\right)  =p^{-1}\int_{0}^{y}\mathbf{F}\left(  z\right)
d\mathbf{G}\left(  z\right)  .$ Since $\mathbf{F}$ and $\mathbf{G}$ are
heavy-tailed, then their right endpoints are infinite and thus they are equal.
Hence, from \cite{W-85}, we may write $\int_{x}^{\infty}d\mathbf{F}\left(
y\right)  /\mathbf{F}\left(  y\right)  =\int_{x}^{\infty}dF\left(  y\right)
/C\left(  y\right)  ,$ where $C\left(  z\right)  :=\mathbf{P}\left(  X\leq
z\leq Y\right)  .$ Differentiating the previous equation leads to the
following crucial equation $C\left(  x\right)  d\mathbf{F}\left(  x\right)
=\mathbf{F}\left(  x\right)  dF\left(  x\right)  ,$ whose solution is defined
by $\mathbf{F}\left(  x\right)  =\exp\left\{  -\int_{x}^{\infty}dF\left(
z\right)  /C\left(  z\right)  \right\}  .$ This leads to Woodroofe's
nonparametric estimator \citep[][]{W-85} of df $\mathbf{F,}$ given by%
\[
\mathbf{F}_{n}^{\left(  \mathbf{W}\right)  }\left(  x\right)  :=\prod
_{i:X_{i}>x}\exp\left\{  -\dfrac{1}{nC_{n}\left(  X_{i}\right)  }\right\}  ,
\]
which is derived only by replacing df's $F$ and $C$ by their respective
empirical counterparts $F_{n}\left(  x\right)  :=n^{-1}\sum\limits_{i=1}%
^{n}\mathbf{1}\left(  X_{i}\leq x\right)  $ and $C_{n}\left(  x\right)
:=n^{-1}\sum\limits_{i=1}^{n}\mathbf{1}\left(  X_{i}\leq x\leq Y_{i}\right)
.$ There exists a more popular estimator for $\mathbf{F,}$ known as
Lynden-Bell nonparametric maximum likelihood estimator \citep[][]{LB-71},
defined by%
\[
\mathbf{F}_{n}^{\left(  \mathbf{LB}\right)  }\left(  x\right)  :=\prod
_{i:X_{i}>x}\left(  1-\frac{1}{nC_{n}\left(  X_{i}\right)  }\right)  ,
\]
which will be considered in this paper to derive a new estimator for the tail
index of df $\mathbf{F}.$ Note that the tail of df $F$ simultaneously depends
on $\overline{\mathbf{G}}$ and $\overline{\mathbf{F}}$ while that of
$\overline{G}$ only relies on $\overline{\mathbf{G}}\mathbf{.}$ By using
Proposition B.1.10 in \cite{deHF06}, to the regularly varying functions
$\overline{\mathbf{F}}$ and $\overline{\mathbf{G}},$ we show that both
$\overline{F}$ and $\overline{G}$ are regularly varying at infinity as well,
with respective indices $-1/\gamma:=-\left(  \gamma_{1}+\gamma_{2}\right)
/\left(  \gamma_{1}\gamma_{2}\right)  $ and $-1/\gamma_{2}.$ In view of the
definition of $\gamma,$ \cite{GS2015} derived a consistent estimator, for the
extreme value index $\gamma_{1},$ whose asymptotic normality is established in
\cite{BchMN-15}, under the tail dependence and the second-order conditions of
regular variation. Recently, by considering a Lynden-Bell integration with a
deterministic threshold $t_{n}>0,$ \cite{WW2016} proposed another
asymptotically normal estimator for $\gamma_{1}$ as follows:%
\[
\widehat{\gamma}_{1}^{\left(  \mathbf{LB}\right)  }\left(  t_{n}\right)
:=\frac{1}{n\overline{\mathbf{F}}_{n}^{\left(  \mathbf{LB}\right)  }\left(
t_{n}\right)  }\sum\limits_{i=1}^{n}\mathbf{1}\left(  X_{i}>t_{n}\right)
\dfrac{\mathbf{F}_{n}^{\left(  \mathbf{LB}\right)  }\left(  X_{i}\right)
}{C_{n}\left(  X_{i}\right)  }\log\dfrac{X_{i}}{t_{n}}.
\]
Likewise, \cite{BchMN-16a} considered a Woodroofe integration (with a random
threshold) to propose a new estimator for the tail index $\gamma_{1}$ given by%
\[
\widehat{\gamma}_{1}^{\left(  \mathbf{W}\right)  }:=\frac{1}{n\overline
{\mathbf{F}}_{n}^{\left(  \mathbf{W}\right)  }\left(  X_{n-k:n}\right)  }%
\sum\limits_{i=1}^{k}\dfrac{\mathbf{F}_{n}^{\left(  \mathbf{W}\right)
}\left(  X_{n-i+1:n}\right)  }{C_{n}\left(  X_{n-i+1:n}\right)  }\log
\dfrac{X_{n-i+1:n}}{X_{n-k:n}},
\]
where, given $n=m=m_{N},$ $Z_{1:m}\leq...\leq Z_{m:m}$ denote the order
statistics pertaining to a sample $Z_{1},...,Z_{m},$ and $k=k_{n}$ is a
(random) sequence of integers such that, given $n=m,$ $1<k_{m}<m,$
$k_{m}\rightarrow\infty$ and $k_{m}/m\rightarrow0$ as $N\rightarrow\infty.$
The consistency and asymptotic normality of $\widehat{\gamma}_{1}^{\left(
\mathbf{W}\right)  }$ are established in \cite{BchMN-16a} through a weak
approximation to Woodroofe's tail process%
\[
\mathbf{D}_{n}^{\left(  \mathbf{W}\right)  }\left(  x\right)  :=\sqrt
{k}\left(  \frac{\overline{\mathbf{F}}_{n}^{\left(  \mathbf{W}\right)
}\left(  X_{n-k:n}x\right)  }{\overline{\mathbf{F}}_{n}^{\left(
\mathbf{W}\right)  }\left(  X_{n-k:n}\right)  }-x^{-1/\gamma_{1}}\right)
,\text{ }x>0.
\]
More precisely, the authors showed that, under $(\ref{second-order})$ and
$(\ref{second-orderG})$ with $\gamma_{1}<\gamma_{2},$ there exist a function
$\mathbf{A}_{0}\left(  t\right)  \sim\mathbf{A}_{\mathbf{F}}^{\ast}\left(
t\right)  :=\mathbf{A}_{\mathbf{F}}\left(  1/\overline{\mathbf{F}}\left(
\mathbb{U}_{F}\left(  t\right)  \right)  \right)  ,$ $t\rightarrow\infty,$ and
a standard Wiener process $\left\{  \mathbf{W}\left(  s\right)  ;\text{ }%
s\geq0\right\}  ,$\ defined on the probability space $\left(  \Omega
,\mathcal{A},\mathbf{P}\right)  ,$\ such that, for $0<\epsilon<1/2-\gamma
/\gamma_{2}$ and $x_{0}>0,$%
\begin{equation}
\sup_{x\geq x_{0}}x^{\left(  1/2-\epsilon\right)  /\gamma-1/\gamma_{2}%
}\left\vert \mathbf{D}_{n}^{\left(  \mathbf{W}\right)  }\left(  x\right)
-\mathbf{\Gamma}\left(  x;\mathbf{W}\right)  -x^{-1/\gamma_{1}}\dfrac
{x^{\tau_{1}/\gamma_{1}}-1}{\gamma_{1}\tau_{1}}\sqrt{k}\mathbf{A}_{0}\left(
n/k\right)  \right\vert =o_{\mathbf{P}}\left(  1\right)  , \label{weak}%
\end{equation}
as $N\rightarrow\infty,$ provided that given $n=m,$ $\sqrt{k_{m}}%
\mathbf{A}_{0}\left(  m/k_{m}\right)  =O\left(  1\right)  ,$\textbf{\ }%
where$\mathbb{\ }\left\{  \Gamma\left(  x;\mathbf{W}\right)  ;\text{
}x>0\right\}  $\textbf{\ }is a Gaussian process defined by%
\begin{align*}
&  \mathbf{\Gamma}\left(  x;\mathbf{W}\right)
\begin{tabular}
[c]{l}%
$:=$%
\end{tabular}
\ \ \ \ \ \ \dfrac{\gamma}{\gamma_{1}}x^{-1/\gamma_{1}}\left\{  x^{1/\gamma
}\mathbf{W}\left(  x^{-1/\gamma}\right)  -\mathbf{W}\left(  1\right)  \right\}
\\
&  \ \ \ \ \ \ \ \ \ \ \ \ \ \ \ \ \ \ +\frac{\gamma}{\gamma_{1}+\gamma_{2}%
}x^{-1/\gamma_{1}}\int_{0}^{1}s^{-\gamma/\gamma_{2}-1}\left\{  x^{1/\gamma
}\mathbf{W}\left(  x^{-1/\gamma}s\right)  -\mathbf{W}\left(  s\right)
\right\}  ds.
\end{align*}
In view of the previous weak approximation, the authors also proved that if,
given $n=m,$ $\sqrt{k_{m}}\mathbf{A}_{\mathbf{F}}^{\ast}\left(  m/k_{m}%
\right)  \rightarrow\lambda,$ then $\sqrt{k}\left(  \widehat{\gamma}%
_{1}-\gamma_{1}\right)  \overset{\mathcal{D}}{\rightarrow}\mathcal{N}\left(
\dfrac{\lambda}{1-\tau_{1}},\sigma^{2}\right)  ,$ as $N\rightarrow\infty,$
where $\sigma^{2}:=\gamma^{2}\left(  1+\gamma_{1}/\gamma_{2}\right)  \left(
1+\left(  \gamma_{1}/\gamma_{2}\right)  ^{2}\right)  /\left(  1-\gamma
_{1}/\gamma_{2}\right)  ^{3}.$ Recently, \cite{BchMN-16b} followed this
approach to introduce a kernel estimator to $\gamma_{1}$ which improves the
bias of $\widehat{\gamma}_{1}^{\left(  \mathbf{W}\right)  }.$ In this paper,
we are interested in Worm's estimator $\widehat{\gamma}_{1}^{\left(
\mathbf{LB}\right)  }\left(  t_{n}\right)  ,$ but with a threshold $t_{n}$
that is assumed to be random and equal to $X_{n-k:n}.$ This makes the
estimator more convenient for numerical implementation than the one with a
deterministic threshold. In other words, we will deal with the following tail
index estimator:%
\[
\widehat{\gamma}_{1}^{\left(  \mathbf{LB}\right)  }:=\frac{1}{n\overline
{\mathbf{F}}_{n}^{\left(  \mathbf{LB}\right)  }\left(  X_{n-k:n}\right)  }%
{\displaystyle\sum\limits_{i=1}^{k}}
\dfrac{\mathbf{F}_{n}^{\left(  \mathbf{LB}\right)  }\left(  X_{n-i+1:n}%
\right)  }{C_{n}\left(  X_{n-i+1:n}\right)  }\log\dfrac{X_{n-i+1:n}}%
{X_{n-k:n}}.
\]
Note that $\mathbf{F}_{n}^{\left(  \mathbf{LB}\right)  }\left(  \infty\right)
=1$ and write $\overline{\mathbf{F}}_{n}^{\left(  \mathbf{LB}\right)  }\left(
X_{n-k:n}\right)  =\int_{X_{n-k:n}}^{\infty}d\mathbf{F}_{n}^{\left(
\mathbf{LB}\right)  }\left(  y\right)  .$ On the other hand, we have
$C_{n}\left(  x\right)  d\mathbf{F}_{n}^{\left(  \mathbf{LB}\right)  }\left(
x\right)  =\mathbf{F}_{n}^{\left(  \mathbf{LB}\right)  }\left(  x\right)
dF_{n}\left(  x\right)  $ (see, e.g., \citeauthor{SS09}, \citeyear[]{SS09}),
then
\[
\overline{\mathbf{F}}_{n}^{\left(  \mathbf{LB}\right)  }\left(  X_{n-k:n}%
\right)  =\int_{X_{n-k:n}}^{\infty}\frac{\mathbf{F}_{n}^{\left(
\mathbf{LB}\right)  }\left(  x\right)  }{C_{n}\left(  x\right)  }dF_{n}\left(
x\right)  =\frac{1}{n}\sum\limits_{i=1}^{k}\frac{\mathbf{F}_{n}^{\left(
\mathbf{LB}\right)  }\left(  X_{n-i+1:n}\right)  }{C_{n}\left(  X_{n-i+1:n}%
\right)  }.
\]
This allows us to rewrite the new estimator into%
\[
\widehat{\gamma}_{1}^{\left(  \mathbf{LB}\right)  }:=%
{\displaystyle\sum\limits_{i=1}^{k}}
a_{n}^{\left(  i\right)  }\log\dfrac{X_{n-i+1:n}}{X_{n-k:n}},
\]
where%
\[
a_{n}^{\left(  i\right)  }:=\frac{\mathbf{F}_{n}^{\left(  \mathbf{LB}\right)
}\left(  X_{n-i+1:n}\right)  }{C_{n}\left(  X_{n-i+1:n}\right)  }/%
{\displaystyle\sum\limits_{i=1}^{k}}
\frac{\mathbf{F}_{n}^{\left(  \mathbf{LB}\right)  }\left(  X_{n-i+1:n}\right)
}{C_{n}\left(  X_{n-i+1:n}\right)  }.
\]
It is worth mentioning that for complete data, we have $n\mathbf{\equiv}N$ and
$\mathbf{F}_{n}\mathbf{\equiv}F_{n}\mathbf{\equiv}C_{n},$ it follows that
$a_{n}^{\left(  i\right)  }\mathbf{\equiv}k^{-1},$ $i=1,...,k$ and
consequently both $\widehat{\gamma}_{1}^{\left(  \mathbf{LB}\right)  }$ and
$\widehat{\gamma}_{1}^{\left(  \mathbf{W}\right)  }$ reduce to the classical
Hill estimator \citep[][]{Hill75}. The consistency and asymptotic normality of
$\widehat{\gamma}_{1}^{\left(  \mathbf{LB}\right)  }$ will be achieved through
a weak approximation of the corresponding tail Lynden-Bell process that we
define by%
\[
\mathbf{D}_{n}^{\left(  \mathbf{LB}\right)  }\left(  x\right)  :=\sqrt
{k}\left(  \frac{\overline{\mathbf{F}}_{n}^{\left(  \mathbf{LB}\right)
}\left(  X_{n-k:n}x\right)  }{\overline{\mathbf{F}}_{n}^{\left(
\mathbf{LB}\right)  }\left(  X_{n-k:n}\right)  }-x^{-1/\gamma_{1}}\right)
,\text{ }x>0.
\]
The rest of the paper is organized as follows. In Section \ref{sec2}, we
provide our main results whose proofs are postponed to Section \ref{sec4}. The
finite sample behavior of the proposed estimator $\widehat{\gamma}%
_{1}^{\left(  \mathbf{LB}\right)  }$ is checked by simulation in Section
\ref{sec3}, where a comparison with the one recently introduced by
\cite{BchMN-16a} is made as well.

\section{\textbf{Main results\label{sec2}}}

\noindent We basically have three main results. The first one, that we give in
Theorem \ref{Theorem1}, consists in an asymptotic relation between the above
mentioned estimators of the distribution tail, namely $\overline{\mathbf{F}%
}_{n}^{\left(  \mathbf{W}\right)  }$ and $\overline{\mathbf{F}}_{n}^{\left(
\mathbf{LB}\right)  }.$ This in turn is instrumental to the Gaussian
approximation of the tail Lynden-Bell process $\mathbf{D}_{n}^{\left(
\mathbf{LB}\right)  }\left(  x\right)  $ stated in Theorem \ref{Theorem2}.
Finally, in Theorem \ref{Theorem3}, we deduce the asymptotic behavior of the
tail index estimator $\widehat{\gamma}_{1}^{\left(  \mathbf{LB}\right)  }.$

\begin{theorem}
\label{Theorem1}Assume that both $\mathbf{F}$ and $\mathbf{G}$ satisfy the
second-order conditions $(\ref{second-order})$ and $(\ref{second-orderG})$
respectively with $\gamma_{1}<\gamma_{2}.$ Let $k=k_{n}$ be a random sequence
of integers such that, given $n=m,$ $k_{m}\rightarrow\infty$ and
$k_{m}/m\rightarrow0,$ as $N\rightarrow\infty,$ then, for any $x_{0}>0,$ we
have%
\[
\sup_{x\geq x_{0}}x^{1/\gamma_{1}}\frac{\left\vert \overline{\mathbf{F}}%
_{n}^{\left(  \mathbf{W}\right)  }\left(  X_{n-k:n}x\right)  -\overline
{\mathbf{F}}_{n}^{\left(  \mathbf{LB}\right)  }\left(  X_{n-k:n}x\right)
\right\vert }{\overline{\mathbf{F}}\left(  X_{n-k:n}\right)  }=O_{\mathbf{P}%
}\left(  \left(  k/n\right)  ^{\gamma_{1}/\gamma}\right)  .
\]

\end{theorem}

\begin{theorem}
\label{Theorem2}Assume that the assumptions of Theorem \ref{Theorem1} hold and
given $n=m,$
\begin{equation}
k_{m}^{1+\gamma_{1}/\left(  2\gamma\right)  }/m\rightarrow0, \label{add}%
\end{equation}
and $\sqrt{k_{m}}\mathbf{A}_{0}\left(  m/k_{m}\right)  =O\left(  1\right)  ,$
as $N\rightarrow\infty.$ Then, for any $x_{0}>0$ and $0<\epsilon
<1/2-\gamma/\gamma_{2},$ we have%
\[
\sup_{x\geq x_{0}}x^{\left(  1/2-\epsilon\right)  /\gamma-1/\gamma_{2}%
}\left\vert \mathbf{D}_{n}^{\left(  \mathbf{LB}\right)  }\left(  x\right)
-\mathbf{\Gamma}\left(  x;\mathbf{W}\right)  -x^{-1/\gamma_{1}}\dfrac
{x^{\tau_{1}/\gamma_{1}}-1}{\gamma_{1}\tau_{1}}\sqrt{k}\mathbf{A}_{0}\left(
n/k\right)  \right\vert =o_{\mathbf{P}}\left(  1\right)  .
\]

\end{theorem}

\begin{theorem}
\label{Theorem3}Assume that $\left(  \ref{RV-1}\right)  $ holds with
$\gamma_{1}<\gamma_{2}$ and let $k=k_{n}$ be a random sequence of integers
such that given $n=m,$ $k_{m}\rightarrow\infty$ and $k_{m}/m\rightarrow0,$ as
$N\rightarrow\infty,$ then $\widehat{\gamma}_{1}^{\left(  \mathbf{LB}\right)
}\overset{\mathbf{P}}{\rightarrow}\gamma_{1}.$ Assume further that the
assumptions of Theorem \ref{Theorem2} hold, then
\begin{align*}
\sqrt{k}\left(  \widehat{\gamma}_{1}^{\left(  \mathbf{LB}\right)  }-\gamma
_{1}\right)   &  =\frac{\sqrt{k}\mathbf{A}_{0}\left(  n/k\right)  }{1-\tau
_{1}}-\gamma\mathbf{W}\left(  1\right) \\
&  +\frac{\gamma}{\gamma_{1}+\gamma_{2}}\int_{0}^{1}\left(  \gamma_{2}%
-\gamma_{1}-\gamma\log s\right)  s^{-\gamma/\gamma_{2}-1}\mathbf{W}\left(
s\right)  ds+o_{\mathbf{P}}\left(  1\right)  .
\end{align*}
If, in addition, we suppose that, given $n=m,$ $\sqrt{k_{m}}\mathbf{A}%
_{\mathbf{F}}^{\ast}\left(  m/k_{m}\right)  \rightarrow\lambda<\infty,$ then
\[
\sqrt{k}\left(  \widehat{\gamma}_{1}^{\left(  \mathbf{LB}\right)  }-\gamma
_{1}\right)  \overset{\mathcal{D}}{\rightarrow}\mathcal{N}\left(
\frac{\lambda}{1-\tau_{1}},\sigma^{2}\right)  ,\text{ as }N\rightarrow\infty.
\]

\end{theorem}

\section{\textbf{Simulation study\label{sec3}}}

\noindent In this section, we illustrate the finite sample behavior of
$\widehat{\gamma}_{1}^{\left(  \mathbf{LB}\right)  }$ and, at the same time,
we compare it with $\widehat{\gamma}_{1}^{\left(  \mathbf{W}\right)  }.$ To
this end, we consider two sets of truncated and truncation data, both drawn
from Burr's model: $\overline{\mathbf{F}}\left(  x\right)  =\left(
1+x^{1/\delta}\right)  ^{-\delta/\gamma_{1}},$ $\overline{\mathbf{G}}\left(
x\right)  =\left(  1+x^{1/\delta}\right)  ^{-\delta/\gamma_{2}},$ $x\geq0,$
where $\delta,\gamma_{1},\gamma_{2}>0.$ The corresponding percentage of
observed data is equal to $p=\gamma_{2}/(\gamma_{1}+\gamma_{2}).$ We fix
$\delta=1/4$ and choose the values $0.6$ and $0.8$ for $\gamma_{1}$ and
$55\%,$ $70\%$ and $90\%$ for $p.$ For each couple $\left(  \gamma
_{1},p\right)  ,$ we solve the equation $p=\gamma_{2}/(\gamma_{1}+\gamma_{2})$
to get the pertaining $\gamma_{2}$-value. We vary the common size $N$ of both
samples $\left(  \mathbf{X}_{1},...,\mathbf{X}_{N}\right)  $ and $\left(
\mathbf{Y}_{1},...,\mathbf{Y}_{N}\right)  ,$ then for each size, we generate
$1000$ independent replicates. Our overall results are taken as the empirical
means of the results obtained through all repetitions. To determine the
optimal number of top statistics used in the computation of the tail index
estimate values, we use the algorithm of \cite{ReTo7}, page 137. Our
illustration and comparison are made with respect to the estimators absolute
biases (abs bias) and the roots of their mean squared errors (rmse). We
summarize the simulation results in Tables \ref{Tab1}, \ref{Tab2} and
\ref{Tab3} for $\gamma_{1}=0.6$ and in Tables \ref{Tab4}, \ref{Tab5} and
\ref{Tab6} for $\gamma_{1}=0.8.$ After the inspection of all the tables, two
conclusions can be drawn regardless of the situation. First, the estimation
accuracy of both estimators decreases when the truncation percentage increases
and this was quite expected. Second, we notice that the newly proposed
estimator $\widehat{\gamma}_{1}^{\left(  \mathbf{LB}\right)  }$ and
$\widehat{\gamma}_{1}^{\left(  \mathbf{W}\right)  }$ behave equally well.%

\begin{table}[tbp] \centering
\begin{tabular}
[c]{cccccccc}%
\multicolumn{8}{c}{$\gamma_{1}=0.6;$ $p=0.55$}\\\hline
&  & \multicolumn{3}{c}{$\widehat{\gamma}_{1}^{\left(  \mathbf{LB}\right)  }$}
& \multicolumn{3}{||c}{$\widehat{\gamma}_{1}^{\left(  \mathbf{W}\right)  }$%
}\\\hline
$N$ & $n$ & abs bias & rmse & $k^{\ast}$ & \multicolumn{1}{||c}{abs bias} &
rmse & $k^{\ast}$\\\hline\hline
\multicolumn{1}{r}{${\small 100}$} & \multicolumn{1}{r}{${\small 54}$} &
\multicolumn{1}{||c}{${\small 0.0407}$} & ${\small 0.2381}$ &
\multicolumn{1}{r}{${\small 26}$} & \multicolumn{1}{||c}{${\small 0.0443}$} &
${\small 0.2328}$ & \multicolumn{1}{r}{${\small 26}$}\\
\multicolumn{1}{r}{${\small 200}$} & \multicolumn{1}{r}{${\small 109}$} &
\multicolumn{1}{||c}{${\small 0.0378}$} & ${\small 0.2610}$ &
\multicolumn{1}{r}{${\small 36}$} & \multicolumn{1}{||c}{${\small 0.0358}$} &
${\small 0.2532}$ & \multicolumn{1}{r}{${\small 37}$}\\
\multicolumn{1}{r}{${\small 300}$} & \multicolumn{1}{r}{${\small 165}$} &
\multicolumn{1}{||c}{${\small 0.0352}$} & ${\small 0.2359}$ &
\multicolumn{1}{r}{${\small 36}$} & \multicolumn{1}{||c}{${\small 0.0323}$} &
${\small 0.2315}$ & \multicolumn{1}{r}{${\small 37}$}\\
\multicolumn{1}{r}{${\small 500}$} & \multicolumn{1}{r}{${\small 274}$} &
\multicolumn{1}{||c}{${\small 0.0199}$} & ${\small 0.2290}$ &
\multicolumn{1}{r}{${\small 61}$} & \multicolumn{1}{||c}{${\small 0.0185}$} &
${\small 0.2238}$ & \multicolumn{1}{r}{${\small 61}$}\\
\multicolumn{1}{r}{${\small 1000}$} & \multicolumn{1}{r}{${\small 549}$} &
\multicolumn{1}{||c}{${\small 0.0074}$} & ${\small 0.1763}$ &
\multicolumn{1}{r}{${\small 112}$} & \multicolumn{1}{||c}{${\small 0.0068}$} &
${\small 0.1748}$ & \multicolumn{1}{r}{${\small 112}$}\\
\multicolumn{1}{r}{${\small 3000}$} & \multicolumn{1}{r}{${\small 1649}$} &
\multicolumn{1}{||c}{${\small 0.0036}$} & ${\small 0.0982}$ &
\multicolumn{1}{r}{${\small 350}$} & \multicolumn{1}{||c}{${\small 0.0037}$} &
${\small 0.0981}$ & \multicolumn{1}{r}{${\small 352}$}\\
\multicolumn{1}{r}{${\small 5000}$} & \multicolumn{1}{r}{${\small 2747}$} &
\multicolumn{1}{||c}{${\small 0.0007}$} & ${\small 0.1066}$ &
\multicolumn{1}{r}{${\small 432}$} & \multicolumn{1}{||c}{${\small 0.0007}$} &
${\small 0.1065}$ & \multicolumn{1}{r}{${\small 432}$}\\\hline\hline
&  &  &  &  &  &  &
\end{tabular}
\caption{Estimation results of Lynden-Bell based (leftt pannel) and Woodroofe based (right pannel) estimators of the shape parameter
$\gamma_{1}=0.6$ of Burr's model through 1000 right-truncated samples with 45$\%$-truncation rate.}\label{Tab1}
\end{table}%
%

\begin{table}[tbp] \centering
\begin{tabular}
[c]{cccccccc}%
\multicolumn{8}{c}{$\gamma_{1}=0.6;$ $p=0.7$}\\\hline
&  & \multicolumn{3}{c}{$\widehat{\gamma}_{1}^{\left(  \mathbf{LB}\right)  }$}
& \multicolumn{3}{||c}{$\widehat{\gamma}_{1}^{\left(  \mathbf{W}\right)  }$%
}\\\hline
$N$ & $n$ & abs bias & rmse & $k^{\ast}$ & \multicolumn{1}{||c}{abs bias} &
rmse & $k^{\ast}$\\\hline\hline
\multicolumn{1}{r}{${\small 100}$} & \multicolumn{1}{r}{${\small 69}$} &
\multicolumn{1}{||c}{${\small 0.0158}$} & ${\small 0.2451}$ &
\multicolumn{1}{r}{${\small 25}$} & \multicolumn{1}{||c}{${\small 0.0144}$} &
${\small 0.2428}$ & \multicolumn{1}{r}{${\small 25}$}\\
\multicolumn{1}{r}{${\small 200}$} & \multicolumn{1}{r}{${\small 140}$} &
\multicolumn{1}{||c}{${\small 0.0095}$} & ${\small 0.1871}$ &
\multicolumn{1}{r}{${\small 39}$} & \multicolumn{1}{||c}{${\small 0.0089}$} &
${\small 0.1866}$ & \multicolumn{1}{r}{${\small 39}$}\\
\multicolumn{1}{r}{${\small 300}$} & \multicolumn{1}{r}{${\small 210}$} &
\multicolumn{1}{||c}{${\small 0.0085}$} & ${\small 0.1590}$ &
\multicolumn{1}{r}{${\small 61}$} & \multicolumn{1}{||c}{${\small 0.0082}$} &
${\small 0.1587}$ & \multicolumn{1}{r}{${\small 61}$}\\
\multicolumn{1}{r}{${\small 500}$} & \multicolumn{1}{r}{${\small 348}$} &
\multicolumn{1}{||c}{${\small 0.0074}$} & ${\small 0.1294}$ &
\multicolumn{1}{r}{${\small 76}$} & \multicolumn{1}{||c}{${\small 0.0072}$} &
${\small 0.1293}$ & \multicolumn{1}{r}{${\small 76}$}\\
\multicolumn{1}{r}{${\small 1000}$} & \multicolumn{1}{r}{${\small 699}$} &
\multicolumn{1}{||c}{${\small 0.0063}$} & ${\small 0.1014}$ &
\multicolumn{1}{r}{${\small 124}$} & \multicolumn{1}{||c}{${\small 0.0062}$} &
${\small 0.1014}$ & \multicolumn{1}{r}{${\small 124}$}\\
\multicolumn{1}{r}{${\small 3000}$} & \multicolumn{1}{r}{${\small 2096}$} &
\multicolumn{1}{||c}{${\small 0.0053}$} & ${\small 0.0962}$ &
\multicolumn{1}{r}{${\small 246}$} & \multicolumn{1}{||c}{${\small 0.0053}$} &
${\small 0.0962}$ & \multicolumn{1}{r}{${\small 246}$}\\
\multicolumn{1}{r}{${\small 5000}$} & \multicolumn{1}{r}{${\small 3498}$} &
\multicolumn{1}{||c}{${\small 0.0036}$} & ${\small 0.0984}$ &
\multicolumn{1}{r}{${\small 400}$} & \multicolumn{1}{||c}{${\small 0.0036}$} &
${\small 0.0984}$ & \multicolumn{1}{r}{${\small 400}$}\\\hline\hline
&  &  &  &  &  &  &
\end{tabular}
\caption{Estimation results of Lynden-Bell based (left pannel) and Woodroofe based (right pannel) estimators of the shape parameter
$\gamma_{1}=0.6$ of Burr's model through 1000 right-truncated samples with 30$\%$-truncation rate.}\label{Tab2}%
\end{table}%
%

\begin{table}[tbp] \centering
\begin{tabular}
[c]{cccccccc}%
\multicolumn{8}{c}{$\gamma_{1}=0.6;$ $p=0.9$}\\\hline
&  & \multicolumn{3}{c}{$\widehat{\gamma}_{1}^{\left(  \mathbf{LB}\right)  }$}
& \multicolumn{3}{||c}{$\widehat{\gamma}_{1}^{\left(  \mathbf{W}\right)  }$%
}\\\hline
$N$ & $n$ & abs bias & rmse & $k^{\ast}$ & \multicolumn{1}{||c}{abs bias} &
rmse & $k^{\ast}$\\\hline\hline
\multicolumn{1}{r}{${\small 100}$} & \multicolumn{1}{r}{${\small 90}$} &
\multicolumn{1}{||c}{${\small 0.0073}$} & ${\small 0.1779}$ &
\multicolumn{1}{r}{${\small 21}$} & \multicolumn{1}{||c}{${\small 0.0070}$} &
${\small 0.1778}$ & \multicolumn{1}{r}{${\small 21}$}\\
\multicolumn{1}{r}{${\small 200}$} & \multicolumn{1}{r}{${\small 180}$} &
\multicolumn{1}{||c}{${\small 0.0066}$} & ${\small 0.1208}$ &
\multicolumn{1}{r}{${\small 54}$} & \multicolumn{1}{||c}{${\small 0.0064}$} &
${\small 0.1208}$ & \multicolumn{1}{r}{${\small 54}$}\\
\multicolumn{1}{r}{${\small 300}$} & \multicolumn{1}{r}{${\small 270}$} &
\multicolumn{1}{||c}{${\small 0.0055}$} & ${\small 0.1133}$ &
\multicolumn{1}{r}{${\small 88}$} & \multicolumn{1}{||c}{${\small 0.0056}$} &
${\small 0.1133}$ & \multicolumn{1}{r}{${\small 88}$}\\
\multicolumn{1}{r}{${\small 500}$} & \multicolumn{1}{r}{${\small 450}$} &
\multicolumn{1}{||c}{${\small 0.0050}$} & ${\small 0.0864}$ &
\multicolumn{1}{r}{${\small 125}$} & \multicolumn{1}{||c}{${\small 0.0050}$} &
${\small 0.0863}$ & \multicolumn{1}{r}{${\small 125}$}\\
\multicolumn{1}{r}{${\small 1000}$} & \multicolumn{1}{r}{${\small 898}$} &
\multicolumn{1}{||c}{${\small 0.0030}$} & ${\small 0.0614}$ &
\multicolumn{1}{r}{${\small 189}$} & \multicolumn{1}{||c}{${\small 0.0029}$} &
${\small 0.0614}$ & \multicolumn{1}{r}{${\small 189}$}\\
\multicolumn{1}{r}{${\small 3000}$} & \multicolumn{1}{r}{${\small 2702}$} &
\multicolumn{1}{||c}{${\small 0.0016}$} & ${\small 0.0494}$ &
\multicolumn{1}{r}{${\small 398}$} & \multicolumn{1}{||c}{${\small 0.0016}$} &
${\small 0.0494}$ & \multicolumn{1}{r}{${\small 398}$}\\
\multicolumn{1}{r}{${\small 5000}$} & \multicolumn{1}{r}{${\small 4496}$} &
\multicolumn{1}{||c}{${\small 0.0010}$} & ${\small 0.0112}$ &
\multicolumn{1}{r}{${\small 467}$} & \multicolumn{1}{||c}{${\small 0.0010}$} &
${\small 0.0112}$ & \multicolumn{1}{r}{${\small 467}$}\\\hline\hline
&  &  &  &  &  &  &
\end{tabular}
\caption{Estimation results of Lynden-Bell based (left pannel) and Woodroofe based (right pannel) estimators of the shape parameter
$\gamma_{1}=0.6$ of Burr's model through 1000 right-truncated samples with 10$\%$-truncation rate.}\label{Tab3}
\end{table}%
%

\begin{table}[tbp] \centering
\begin{tabular}
[c]{cccccccc}%
\multicolumn{8}{c}{$\gamma_{1}=0.8;$ $p=0.55$}\\\hline
&  & \multicolumn{3}{c}{$\widehat{\gamma}_{1}^{\left(  \mathbf{LB}\right)  }$}
& \multicolumn{3}{||c}{$\widehat{\gamma}_{1}^{\left(  \mathbf{W}\right)  }$%
}\\\hline
$N$ & $n$ & abs bias & rmse & $k^{\ast}$ & \multicolumn{1}{||c}{abs bias} &
rmse & $k^{\ast}$\\\hline\hline
\multicolumn{1}{r}{${\small 100}$} & \multicolumn{1}{r}{${\small 55}$} &
\multicolumn{1}{||c}{${\small 0.0570}$} & ${\small 0.3330}$ &
\multicolumn{1}{r}{${\small 30}$} & \multicolumn{1}{||c}{${\small 0.0636}$} &
${\small 0.3167}$ & \multicolumn{1}{r}{${\small 31}$}\\
\multicolumn{1}{r}{${\small 200}$} & \multicolumn{1}{r}{${\small 110}$} &
\multicolumn{1}{||c}{${\small 0.0401}$} & ${\small 0.3604}$ &
\multicolumn{1}{r}{${\small 33}$} & \multicolumn{1}{||c}{${\small 0.0347}$} &
${\small 0.3453}$ & \multicolumn{1}{r}{${\small 35}$}\\
\multicolumn{1}{r}{${\small 300}$} & \multicolumn{1}{r}{${\small 164}$} &
\multicolumn{1}{||c}{${\small 0.0252}$} & ${\small 0.2563}$ &
\multicolumn{1}{r}{${\small 69}$} & \multicolumn{1}{||c}{${\small 0.0272}$} &
${\small 0.2530}$ & \multicolumn{1}{r}{${\small 71}$}\\
\multicolumn{1}{r}{${\small 500}$} & \multicolumn{1}{r}{${\small 276}$} &
\multicolumn{1}{||c}{${\small 0.0227}$} & ${\small 0.1807}$ &
\multicolumn{1}{r}{${\small 112}$} & \multicolumn{1}{||c}{${\small 0.0216}$} &
${\small 0.1794}$ & \multicolumn{1}{r}{${\small 113}$}\\
\multicolumn{1}{r}{${\small 1000}$} & \multicolumn{1}{r}{${\small 551}$} &
\multicolumn{1}{||c}{${\small 0.0148}$} & ${\small 0.1795}$ &
\multicolumn{1}{r}{${\small 196}$} & \multicolumn{1}{||c}{${\small 0.0142}$} &
${\small 0.1788}$ & \multicolumn{1}{r}{${\small 197}$}\\
\multicolumn{1}{r}{${\small 3000}$} & \multicolumn{1}{r}{${\small 1647}$} &
\multicolumn{1}{||c}{${\small 0.0124}$} & ${\small 0.1794}$ &
\multicolumn{1}{r}{${\small 525}$} & \multicolumn{1}{||c}{${\small 0.0121}$} &
${\small 0.1783}$ & \multicolumn{1}{r}{${\small 525}$}\\
\multicolumn{1}{r}{${\small 5000}$} & \multicolumn{1}{r}{${\small 2751}$} &
\multicolumn{1}{||c}{${\small 0.0075}$} & ${\small 0.1260}$ &
\multicolumn{1}{r}{${\small 688}$} & \multicolumn{1}{||c}{${\small 0.0074}$} &
${\small 0.1259}$ & \multicolumn{1}{r}{${\small 688}$}\\\hline\hline
&  &  &  &  &  &  &
\end{tabular}
\caption{Estimation results of Lynden-Bell based (left pannel) and Woodroofe based (right pannel) estimators of the shape parameter
$\gamma_{1}=0.8$ of Burr's model through 1000 right-truncated samples with 45$\%$-truncation rate.}\label{Tab4}
\end{table}%
%

\begin{table}[tbp] \centering
\begin{tabular}
[c]{cccccccc}%
\multicolumn{8}{c}{$\gamma_{1}=0.8;$ $p=0.7$}\\\hline
&  & \multicolumn{3}{c}{$\widehat{\gamma}_{1}^{\left(  \mathbf{LB}\right)  }$}
& \multicolumn{3}{||c}{$\widehat{\gamma}_{1}^{\left(  \mathbf{W}\right)  }$%
}\\\hline
$N$ & $n$ & abs bias & rmse & $k^{\ast}$ & \multicolumn{1}{||c}{abs bias} &
rmse & $k^{\ast}$\\\hline\hline
\multicolumn{1}{r}{${\small 100}$} & \multicolumn{1}{r}{${\small 69}$} &
\multicolumn{1}{||c}{${\small 0.0217}$} & ${\small 0.3827}$ &
\multicolumn{1}{r}{${\small 28}$} & \multicolumn{1}{||c}{${\small 0.0195}$} &
${\small 0.3787}$ & \multicolumn{1}{r}{${\small 28}$}\\
\multicolumn{1}{r}{${\small 200}$} & \multicolumn{1}{r}{${\small 139}$} &
\multicolumn{1}{||c}{${\small 0.0203}$} & ${\small 0.2918}$ &
\multicolumn{1}{r}{${\small 59}$} & \multicolumn{1}{||c}{${\small 0.0194}$} &
${\small 0.2905}$ & \multicolumn{1}{r}{${\small 59}$}\\
\multicolumn{1}{r}{${\small 300}$} & \multicolumn{1}{r}{${\small 210}$} &
\multicolumn{1}{||c}{${\small 0.0189}$} & ${\small 0.1857}$ &
\multicolumn{1}{r}{${\small 66}$} & \multicolumn{1}{||c}{${\small 0.0184}$} &
${\small 0.1852}$ & \multicolumn{1}{r}{${\small 66}$}\\
\multicolumn{1}{r}{${\small 500}$} & \multicolumn{1}{r}{${\small 348}$} &
\multicolumn{1}{||c}{${\small 0.0143}$} & ${\small 0.1593}$ &
\multicolumn{1}{r}{${\small 113}$} & \multicolumn{1}{||c}{${\small 0.0140}$} &
${\small 0.1591}$ & \multicolumn{1}{r}{${\small 113}$}\\
\multicolumn{1}{r}{${\small 1000}$} & \multicolumn{1}{r}{${\small 700}$} &
\multicolumn{1}{||c}{${\small 0.0049}$} & ${\small 0.1205}$ &
\multicolumn{1}{r}{${\small 230}$} & \multicolumn{1}{||c}{${\small 0.0049}$} &
${\small 0.1204}$ & \multicolumn{1}{r}{${\small 230}$}\\
\multicolumn{1}{r}{${\small 3000}$} & \multicolumn{1}{r}{${\small 2100}$} &
\multicolumn{1}{||c}{${\small 0.0037}$} & ${\small 0.0886}$ &
\multicolumn{1}{r}{${\small 449}$} & \multicolumn{1}{||c}{${\small 0.0038}$} &
${\small 0.0886}$ & \multicolumn{1}{r}{${\small 449}$}\\
\multicolumn{1}{r}{${\small 5000}$} & \multicolumn{1}{r}{${\small 3500}$} &
\multicolumn{1}{||c}{${\small 0.0031}$} & ${\small 0.0857}$ &
\multicolumn{1}{r}{${\small 500}$} & \multicolumn{1}{||c}{${\small 0.0031}$} &
${\small 0.0857}$ & \multicolumn{1}{r}{${\small 500}$}\\\hline\hline
&  &  &  &  &  &  &
\end{tabular}
\caption{Estimation results of Lynden-Bell based (left pannel) and Woodroofe based (right pannel) estimators of the shape parameter
$\gamma_{1}=0.8$ of Burr's model through 1000 right-truncated samples with 30$\%$-truncation rate.}\label{Tab5}
\end{table}%
%

\begin{table}[tbp] \centering
\begin{tabular}
[c]{cccccccc}%
\multicolumn{8}{c}{$\gamma_{1}=0.8;$ $p=0.9$}\\\hline
&  & \multicolumn{3}{c}{$\widehat{\gamma}_{1}^{\left(  \mathbf{LB}\right)  }$}
& \multicolumn{3}{||c}{$\widehat{\gamma}_{1}^{\left(  \mathbf{W}\right)  }$%
}\\\hline
$N$ & $n$ & abs bias & rmse & $k^{\ast}$ & \multicolumn{1}{||c}{abs bias} &
rmse & $k^{\ast}$\\\hline\hline
\multicolumn{1}{r}{${\small 100}$} & \multicolumn{1}{r}{${\small 89}$} &
\multicolumn{1}{||c}{${\small 0.0380}$} & ${\small 0.1833}$ &
\multicolumn{1}{r}{${\small 38}$} & \multicolumn{1}{||c}{${\small 0.0369}$} &
${\small 0.1827}$ & \multicolumn{1}{r}{${\small 38}$}\\
\multicolumn{1}{r}{${\small 200}$} & \multicolumn{1}{r}{${\small 179}$} &
\multicolumn{1}{||c}{${\small 0.0345}$} & ${\small 0.1383}$ &
\multicolumn{1}{r}{${\small 80}$} & \multicolumn{1}{||c}{${\small 0.0342}$} &
${\small 0.1383}$ & \multicolumn{1}{r}{${\small 80}$}\\
\multicolumn{1}{r}{${\small 300}$} & \multicolumn{1}{r}{${\small 269}$} &
\multicolumn{1}{||c}{${\small 0.0173}$} & ${\small 0.1014}$ &
\multicolumn{1}{r}{${\small 99}$} & \multicolumn{1}{||c}{${\small 0.0175}$} &
${\small 0.1013}$ & \multicolumn{1}{r}{${\small 99}$}\\
\multicolumn{1}{r}{${\small 500}$} & \multicolumn{1}{r}{${\small 450}$} &
\multicolumn{1}{||c}{${\small 0.0108}$} & ${\small 0.0927}$ &
\multicolumn{1}{r}{${\small 143}$} & \multicolumn{1}{||c}{${\small 0.0106}$} &
${\small 0.0926}$ & \multicolumn{1}{r}{${\small 143}$}\\
\multicolumn{1}{r}{${\small 1000}$} & \multicolumn{1}{r}{${\small 899}$} &
\multicolumn{1}{||c}{${\small 0.0021}$} & ${\small 0.0729}$ &
\multicolumn{1}{r}{${\small 260}$} & \multicolumn{1}{||c}{${\small 0.0021}$} &
${\small 0.0729}$ & \multicolumn{1}{r}{${\small 260}$}\\
\multicolumn{1}{r}{${\small 3000}$} & \multicolumn{1}{r}{${\small 2697}$} &
\multicolumn{1}{||c}{${\small 0.0013}$} & ${\small 0.0591}$ &
\multicolumn{1}{r}{${\small 443}$} & \multicolumn{1}{||c}{${\small 0.0013}$} &
${\small 0.0591}$ & \multicolumn{1}{r}{${\small 443}$}\\
\multicolumn{1}{r}{${\small 5000}$} & \multicolumn{1}{r}{${\small 4500}$} &
\multicolumn{1}{||c}{${\small 0.0001}$} & ${\small 0.0309}$ &
\multicolumn{1}{r}{${\small 997}$} & \multicolumn{1}{||c}{${\small 0.0001}$} &
${\small 0.0309}$ & \multicolumn{1}{r}{${\small 997}$}\\\hline\hline
&  &  &  &  &  &  &
\end{tabular}
\caption{Estimation results of Lynden-Bell based (left pannel) and Woodroofe based (right pannel) estimators of the shape parameter
$\gamma_{1}=0.8$ of Burr's model through 1000 right-truncated samples with 10$\%$-truncation rate.}\label{Tab6}
\end{table}%
\bigskip

\section{\textbf{Proofs\label{sec4}}}

\subsection{Proof Theorem \ref{Theorem1}}

For $x\geq x_{0}$ we have%
\[
\mathbf{F}_{n}^{\left(  \mathbf{W}\right)  }\left(  X_{n-k:n}x\right)
=\exp\left\{  -\int_{X_{n-k:n}x}^{\infty}\frac{dF_{n}\left(  y\right)  }%
{C_{n}\left(  y\right)  }\right\}  .
\]
We show that the latter exponent is negligible in probability uniformly over
$x\geq x_{0}.$ Indeed, note that both $F_{n}\left(  y\right)  /F\left(
y\right)  $ and $C\left(  y\right)  /C_{n}\left(  y\right)  $ are
stochastically bounded from above on $y<X_{n:n}$ (see, e.g.,
\citeauthor{Shorack-86}, \citeyear[page
415]{Shorack-86} and \citeauthor{SS09}, \citeyear[]{SS09}, respectively), it
follows that%
\begin{equation}
-\int_{X_{n-k:n}x}^{\infty}\frac{dF_{n}\left(  y\right)  }{C_{n}\left(
y\right)  }=O_{\mathbf{P}}\left(  1\right)  \int_{X_{n-k:n}x}^{\infty}%
\frac{dF\left(  y\right)  }{C\left(  y\right)  }. \label{equ0}%
\end{equation}
By a change of variables we have%
\begin{equation}
\int_{X_{n-k:n}x}^{\infty}\frac{d\overline{F}\left(  y\right)  }{C\left(
y\right)  }=\frac{\overline{F}\left(  X_{n-k:n}\right)  }{C\left(
X_{n-k:n}\right)  }\left(  \int_{x}^{\infty}\frac{C\left(  X_{n-k:n}\right)
}{C\left(  X_{n-k:n}t\right)  }d\frac{\overline{F}\left(  X_{n-k:n}t\right)
}{\overline{F}\left(  X_{n-k:n}\right)  }\right)  . \label{equ1}%
\end{equation}
Recall that $X_{n-k:n}\overset{\mathbf{P}}{\rightarrow}\infty$ and that
$\overline{F}$ is regularly varying at infinity with index $-1/\gamma.$ On the
other hand, from Assertion $\left(  i\right)  $ of Lemma A.2 in
\cite{BchMN-16a} we deduce that $1/C$ is also regularly varying at infinity
with index $1/\gamma_{2}.$ Thus, we may apply Potters inequalities, see e.g.
Proposition B.1.10 in \cite{deHF06}, to both $\overline{F}$ and $1/C$ to
write: for all large $N,$ any $t\geq x_{0}$ and any sufficiently small
$\delta,\nu>0,$ with large probability,
\begin{equation}
\left\vert \frac{\overline{F}\left(  X_{n-k:n}t\right)  }{\overline{F}\left(
X_{n-k:n}\right)  }-t^{-1/\gamma}\right\vert <\delta t^{-1/\gamma\pm\nu}\text{
and }\left\vert \frac{C\left(  X_{n-k:n}\right)  }{C\left(  X_{n-k:n}t\right)
}-t^{1/\gamma_{2}}\right\vert <\delta t^{1/\gamma_{2}\pm\nu}, \label{pot}%
\end{equation}
where $t^{\pm a}:=\max\left(  t^{a},t^{-a}\right)  .$ These two inequalities
may be rewritten, into%
\[
\frac{\overline{F}\left(  X_{n-k:n}t\right)  }{\overline{F}\left(
X_{n-k:n}\right)  }=t^{-1/\gamma}\left(  1+o_{\mathbf{P}}\left(  t^{\pm\nu
}\right)  \right)  \text{ and }\frac{C\left(  X_{n-k:n}\right)  }{C\left(
X_{n-k:n}t\right)  }=t^{1/\gamma_{2}}\left(  1+o_{\mathbf{P}}\left(  t^{\pm
\nu}\right)  \right)  ,
\]
uniformly on $t\geq x_{0}.$ This leads to%
\begin{equation}
\int_{x}^{\infty}\frac{C\left(  X_{n-k:n}\right)  }{C\left(  X_{n-k:n}%
t\right)  }d\frac{\overline{F}\left(  X_{n-k:n}t\right)  }{\overline{F}\left(
X_{n-k:n}\right)  }=-\frac{\gamma_{1}}{\gamma}x^{-1/\gamma_{1}}\left(
1+o_{\mathbf{P}}\left(  x^{\pm\nu}\right)  \right)  . \label{equ2}%
\end{equation}
In view of $(\ref{HJ}),$ \cite{BchMN-16a} showed, in Lemma A1, that
$\overline{F}\left(  y\right)  =\left(  1+o\left(  1\right)  \right)
c_{1}y^{-1/\gamma}$ and $\overline{G}\left(  y\right)  =\left(  1+o\left(
1\right)  \right)  c_{2}y^{-1/\gamma_{2}}$ as $y\rightarrow\infty,$ for some
constants $c_{1},c_{2}>0.$ In other words, $\mathbb{U}_{F}\left(  s\right)
=\left(  1+o\left(  1\right)  \right)  \left(  c_{1}s\right)  ^{\gamma}$as
$s\rightarrow\infty$, and $C\left(  y\right)  =\left(  1+o\left(  1\right)
\right)  c_{2}y^{-1/\gamma_{2}}$ as $y\rightarrow\infty.$ On the other hand,
from Lemma A4 in \cite{BchMN-16a}, we have $X_{n-k:n}=\left(  1+o_{\mathbf{P}%
}\left(  1\right)  \right)  \mathbb{U}_{F}\left(  n/k\right)  ,$ it follows
that $X_{n-k:n}=\left(  1+o_{\mathbf{P}}\left(  1\right)  \right)
c_{1}^{\gamma}\left(  k/n\right)  ^{-\gamma}.$ Note that $1-\gamma/\gamma
_{2}=\gamma/\gamma_{1},$ hence%
\begin{equation}
\frac{\overline{F}\left(  X_{n-k:n}\right)  }{C\left(  X_{n-k:n}\right)
}=\left(  1+o_{\mathbf{P}}\left(  1\right)  \right)  c_{1}^{\gamma/\gamma_{2}%
}c_{2}^{-1}\left(  k/n\right)  ^{\gamma/\gamma_{1}}. \label{equ3}%
\end{equation}
Plugging results $\left(  \ref{equ2}\right)  $ and $\left(  \ref{equ3}\right)
$ in equation $\left(  \ref{equ1}\right)  $ yields%
\begin{equation}
\int_{X_{n-k:n}x}^{\infty}\frac{d\overline{F}\left(  y\right)  }{C\left(
y\right)  }=\left(  k/n\right)  ^{\gamma/\gamma_{1}}c_{1}^{\gamma/\gamma_{2}%
}c_{2}^{-1}\gamma_{1}x^{-1/\gamma_{1}}\left(  1+o_{\mathbf{P}}\left(
x^{\pm\nu}\right)  \right)  . \label{int-o}%
\end{equation}
By combining equations $(\ref{equ0})$ and $(\ref{int-o}),$ we obtain%
\begin{equation}
\int_{X_{n-k:n}x}^{\infty}\frac{dF_{n}\left(  y\right)  }{C_{n}\left(
y\right)  }=O_{\mathbf{P}}\left(  1\right)  \left(  k/n\right)  ^{\gamma
/\gamma_{1}}x^{-1/\gamma_{1}}\left(  1+o_{\mathbf{P}}\left(  x^{\pm\nu
}\right)  \right)  , \label{int}%
\end{equation}
which obviously tends to zero in probability (uniformly on $x\geq x_{0}).$ We
may now apply Taylor's expansion $e^{t}=1+t+O\left(  t^{2}\right)  ,$ as
$t\rightarrow0,$ to get%
\[
\exp\left\{  -\int_{X_{n-k:n}x}^{\infty}\frac{dF_{n}\left(  y\right)  }%
{C_{n}\left(  y\right)  }\right\}  =1-\int_{X_{n-k:n}x}^{\infty}\frac
{dF_{n}\left(  y\right)  }{C_{n}\left(  y\right)  }+O_{\mathbf{P}}\left(
\int_{X_{n-k:n}x}^{\infty}\frac{dF_{n}\left(  y\right)  }{C_{n}\left(
y\right)  }\right)  ^{2},\text{ }N\rightarrow\infty.
\]
In other words, we have%
\begin{equation}
\overline{\mathbf{F}}_{n}^{\left(  \mathbf{W}\right)  }\left(  X_{n-k:n}%
x\right)  =\int_{X_{n-k:n}x}^{\infty}\frac{dF_{n}\left(  y\right)  }%
{C_{n}\left(  y\right)  }+R_{n1}\left(  x\right)  ,\text{ }N\rightarrow\infty,
\label{11}%
\end{equation}
where $R_{n1}\left(  x\right)  :=O_{\mathbf{P}}\left(  \left(  k/n\right)
^{2\gamma/\gamma_{1}}\right)  x^{-2/\gamma_{1}}\left(  1+o_{\mathbf{P}}\left(
x^{\pm\nu}\right)  \right)  .$ Next, we show that%
\begin{equation}
\overline{\mathbf{F}}_{n}^{\left(  \mathbf{LB}\right)  }\left(  X_{n-k:n}%
x\right)  =\int_{X_{n-k:n}x}^{\infty}\frac{dF_{n}\left(  y\right)  }%
{C_{n}\left(  y\right)  }+R_{n2}\left(  x\right)  ,\text{ }N\rightarrow\infty.
\label{12}%
\end{equation}
Observe that, by taking the logarithme then its exponential in the definition
of $\mathbf{F}_{n}^{\left(  \mathbf{LB}\right)  }\left(  x\right)  ,$ we have%
\[
\mathbf{F}_{n}^{\left(  \mathbf{LB}\right)  }\left(  X_{n-k:n}x\right)
=\exp\left\{  \sum_{i=1}^{n}\mathbf{1}\left(  X_{i:n}>X_{n-k:n}x\right)
\log\left(  1-\dfrac{1}{nC_{n}\left(  X_{i:n}\right)  }\right)  \right\}  ,
\]
which may be rewritten into $\exp\left\{  n\int_{x}^{\infty}\log\left(
1-\dfrac{1}{nC_{n}\left(  X_{n-k:n}y\right)  }\right)  dF_{n}\left(
X_{n-k:n}y\right)  \right\}  .$ To get approximation $(\ref{12})$ it suffices
to apply successively, in the previous quantity, Taylor's expansions
$e^{t}=1+t+O\left(  t^{2}\right)  $ and $\log\left(  1-t\right)  =-t+O\left(
t^{2}\right)  $ (as $t\rightarrow0)$ with similar arguments as above (we omit
further details). Combining $(\ref{11})$ and $(\ref{12})$ and setting
$R_{n}\left(  x\right)  :=R_{n1}\left(  x\right)  -R_{n2}\left(  x\right)  $
yield%
\begin{equation}
\overline{\mathbf{F}}_{n}^{\left(  \mathbf{W}\right)  }\left(  X_{n-k:n}%
x\right)  -\overline{\mathbf{F}}_{n}^{\left(  \mathbf{LB}\right)  }\left(
X_{n-k:n}x\right)  =R_{n}\left(  x\right)  ,\text{ }N\rightarrow\infty.
\label{diff}%
\end{equation}
On the other hand, by once again using Taylor's expansion, we write%
\[
\overline{\mathbf{F}}\left(  X_{n-k:n}\right)  =\int_{X_{n-k:n}}^{\infty}%
\frac{dF\left(  y\right)  }{C\left(  y\right)  }+\widetilde{R}_{n}\left(
x\right)  ,\text{ }N\rightarrow\infty.
\]
From equation $(\ref{int-o})$, we infer that $\overline{\mathbf{F}}\left(
X_{n-k:n}\right)  =c_{2}^{-1}c_{1}^{1-\gamma/\gamma_{1}}\left(  k/n\right)
^{\gamma/\gamma_{1}}\left(  1+o_{\mathbf{P}}\left(  1\right)  \right)  ,$
which implies, in view of $(\ref{diff}),$ that%
\[
x^{1/\gamma_{1}}\frac{\overline{\mathbf{F}}_{n}^{\left(  \mathbf{LB}\right)
}\left(  X_{n-k:n}x\right)  -\overline{\mathbf{F}}_{n}^{\left(  \mathbf{W}%
\right)  }\left(  X_{n-k:n}x\right)  }{\overline{\mathbf{F}}\left(
X_{n-k:n}\right)  }=O_{\mathbf{P}}\left(  \left(  k/n\right)  ^{\gamma
/\gamma_{1}}\right)  x^{-1/\gamma_{1}\pm\nu}.
\]
Observe now that, for a sufficiently small $\nu>0,$ we have $x^{-1/\gamma
_{1}\pm\nu}=O_{\mathbf{P}}\left(  1\right)  ,$ uniformly on $x\geq x_{0}>0,$
as sought.

\subsection{Proof Theorem \ref{Theorem2}}

In a similar way to what is done with $\mathbf{D}_{n}^{\left(  \mathbf{W}%
\right)  }\left(  x\right)  ,$ in the proof of Theorem 2.1 in \cite{BchMN-16a}%
, we decompose $k^{-1/2}\mathbf{D}_{n}^{\left(  \mathbf{LB}\right)  }\left(
x\right)  $ into the sum of%
\[
\mathbf{N}_{n1}\left(  x\right)  :=x^{-1/\gamma_{1}}\frac{\overline
{\mathbf{F}}_{n}^{\left(  \mathbf{LB}\right)  }\left(  X_{n-k:n}x\right)
-\overline{\mathbf{F}}\left(  X_{n-k:n}\right)  }{\overline{\mathbf{F}}\left(
X_{n-k:n}\right)  },
\]%
\[
\mathbf{N}_{n2}\left(  x\right)  :=-\frac{\overline{\mathbf{F}}\left(
X_{n-k:n}x\right)  }{\overline{\mathbf{F}}_{n}^{\left(  \mathbf{LB}\right)
}\left(  X_{n-k:n}\right)  }\frac{\overline{\mathbf{F}}_{n}^{\left(
\mathbf{LB}\right)  }\left(  X_{n-k:n}\right)  -\overline{\mathbf{F}}\left(
X_{n-k:n}\right)  }{\overline{\mathbf{F}}\left(  X_{n-k:n}\right)  },\medskip
\]%
\[
\mathbf{N}_{n3}\left(  x\right)  :=\left(  \frac{\overline{\mathbf{F}}\left(
X_{n-k:n}x\right)  }{\overline{\mathbf{F}}_{n}^{\left(  \mathbf{LB}\right)
}\left(  X_{n-k:n}\right)  }-x^{-1/\gamma_{1}}\right)  \frac{\overline
{\mathbf{F}}_{n}^{\left(  \mathbf{LB}\right)  }\left(  X_{n-k:n}x\right)
-\overline{\mathbf{F}}\left(  xX_{n-k:n}\right)  }{\overline{\mathbf{F}%
}\left(  X_{n-k:n}x\right)  },\medskip
\]
and $\mathbf{N}_{n4}\left(  x\right)  :=\overline{\mathbf{F}}\left(
X_{n-k:n}x\right)  /\overline{\mathbf{F}}\left(  X_{n-k:n}\right)
-x^{-1/\gamma_{1}}.$ If we let%
\[
\mathbf{M}_{n1}\left(  x\right)  :=x^{-1/\gamma_{1}}\frac{\overline
{\mathbf{F}}_{n}^{\left(  \mathbf{W}\right)  }\left(  X_{n-k:n}x\right)
-\overline{\mathbf{F}}\left(  X_{n-k:n}\right)  }{\overline{\mathbf{F}}\left(
X_{n-k:n}\right)  },
\]
then, by applying Theorem \ref{Theorem1}, we have $x^{1/\gamma_{1}}%
\mathbf{N}_{n1}\left(  x\right)  =x^{1/\gamma_{1}}\mathbf{M}_{n1}\left(
x\right)  +x^{-1/\gamma_{1}}O_{\mathbf{P}}\left(  \left(  k/n\right)
^{\gamma/\gamma_{1}}\right)  ,$ uniformly on $x\geq x_{0}.$ By assumption we
have $k^{1+\gamma_{1}/\left(  2\gamma\right)  }/n\overset{\mathbf{P}%
}{\rightarrow}0,$ which is equivalent to $\sqrt{k}\left(  k/n\right)
^{\gamma/\gamma_{1}}\overset{\mathbf{P}}{\rightarrow}0$ as $N\rightarrow
\infty,$ therefore%
\begin{equation}
x^{1/\gamma_{1}}\sqrt{k}\mathbf{N}_{n1}\left(  x\right)  =x^{1/\gamma_{1}%
}\sqrt{k}\mathbf{M}_{n1}\left(  x\right)  +o_{\mathbf{P}}\left(
x^{-1/\gamma_{1}}\right)  . \label{rep}%
\end{equation}
In view of this representation we show that, both $\mathbf{D}_{n}^{\left(
\mathbf{W}\right)  }\left(  x\right)  $ and $\mathbf{D}_{n}^{\left(
\mathbf{LB}\right)  }\left(  x\right)  $ are (weakly) approximated, in the
probability space $\left(  \Omega,\mathcal{A},\mathbf{P}\right)  ,$ by the
same Gaussian process $\mathbf{\Gamma}\left(  x;\mathbf{W}\right)  $ given in
$\left(  \ref{weak}\right)  .$ Indeed, for a sufficiently small $\epsilon>0,$
and $0<\eta<1/2,$ \cite{BchMN-16a} (see the beginning of the proof of Theorem
2.1 therein), showed that
\[
x^{1/\gamma_{1}}\sqrt{k}\mathbf{M}_{n1}\left(  x\right)  =\Phi\left(
x\right)  +o_{\mathbf{P}}\left(  x^{\left(  1-\eta\right)  /\gamma\pm\epsilon
}\right)  ,
\]
where $\Phi\left(  x\right)  :=x^{1/\gamma}\left\{  \dfrac{\gamma}{\gamma_{1}%
}\mathbf{W}\left(  x^{-1/\gamma}\right)  +\dfrac{\gamma}{\gamma_{1}+\gamma
_{2}}\int_{0}^{1}t^{-\gamma/\gamma_{2}-1}\mathbf{W}\left(  x^{-1/\gamma
}t\right)  dt\right\}  .$ Then by using representation $\left(  \ref{rep}%
\right)  ,$ we get $x^{1/\gamma_{1}}\sqrt{k}\mathbf{N}_{n1}\left(  x\right)
=\Phi\left(  x\right)  +o_{\mathbf{P}}\left(  x^{-1/\gamma_{1}}\right)
+o_{\mathbf{P}}\left(  x^{\left(  1-\eta\right)  /\gamma\pm\epsilon}\right)
.$ In particular for $x=1,$ we have%
\begin{equation}
\sqrt{k}\left(  \frac{\overline{\mathbf{F}}_{n}^{\left(  \mathbf{LB}\right)
}\left(  X_{n-k:n}\right)  }{\overline{\mathbf{F}}\left(  X_{n-k:n}\right)
}-1\right)  =\sqrt{k}\mathbf{N}_{n1}\left(  1\right)  =\Phi\left(  1\right)
+o_{\mathbf{P}}\left(  1\right)  , \label{N1}%
\end{equation}
leading to $\overline{\mathbf{F}}_{n}^{\left(  \mathbf{LB}\right)  }\left(
X_{n-k:n}\right)  /\overline{\mathbf{F}}\left(  X_{n-k:n}\right)
\overset{\mathbf{P}}{\rightarrow}1,$ as $N\rightarrow\infty.$ By applying
Potters inequalities to $\overline{\mathbf{F}}$ (as it was done for
$\overline{F}$ in $(\ref{equ3}))$ together with the previous limit, we obtain
\begin{equation}
\frac{\overline{\mathbf{F}}\left(  X_{n-k:n}x\right)  }{\overline{\mathbf{F}%
}_{n}^{\left(  \mathbf{LB}\right)  }\left(  X_{n-k:n}\right)  }=\left(
1+O_{\mathbf{P}}\left(  x^{\pm\epsilon}\right)  \right)  x^{-1/\gamma_{1}}.
\label{rap}%
\end{equation}
By combining $\left(  \ref{N1}\right)  $ and$\left(  \ref{rap}\right)  ,$ we
get $x^{1/\gamma_{1}}\sqrt{k}\mathbf{N}_{n2}\left(  x\right)  =-\Phi\left(
1\right)  +o_{\mathbf{P}}\left(  x^{\pm\epsilon}\right)  .$ For the third term
$\mathbf{N}_{n3}\left(  x\right)  ,$ we use similar arguments to show that
\[
x^{1/\gamma_{1}}\sqrt{k}\mathbf{N}_{n3}\left(  x\right)  =o_{\mathbf{P}%
}\left(  x^{-1/\gamma_{1}\pm\epsilon}\right)  +o_{\mathbf{P}}\left(
x^{-1/\gamma_{1}+\left(  1-\eta\right)  /\gamma\pm\epsilon}\right)  .
\]
Observe that $x^{1/\gamma_{1}-\left(  1-\eta_{0}\right)  /\gamma}%
o_{\mathbf{P}}\left(  x^{-1/\gamma_{1}\pm\epsilon}\right)  $ and
$x^{1/\gamma_{1}-\left(  1-\eta_{0}\right)  /\gamma}o_{\mathbf{P}}\left(
x^{-1/\gamma_{1}+\left(  1-\eta\right)  /\gamma\pm\epsilon}\right)  $
respectively equal $o_{\mathbf{P}}\left(  x^{-\left(  1-\eta_{0}\right)
/\gamma\pm\epsilon}\right)  $ and $o_{\mathbf{P}}\left(  x^{\left(  \eta
-\eta_{0}\right)  /\gamma\pm\epsilon}\right)  ,$ for $\gamma/\gamma_{2}%
<\eta_{0}<\eta<1/2,$ and that both the last two quantities are equal to
$o_{\mathbf{P}}\left(  1\right)  $ for any small $\epsilon>0$ and $x\geq
x_{0}>0.$ Finally, by following the same steps at the end of the proof of
Theorem 2.1 in \cite{BchMN-16a}, we get%

\[
\sqrt{k}\mathbf{N}_{n4}\left(  x\right)  =x^{-1/\gamma_{1}}\dfrac{x^{\tau
_{1}/\gamma_{1}}-1}{\gamma_{1}\tau_{1}}\sqrt{k}\mathbf{A}_{0}\left(
n/k\right)  +o_{\mathbf{p}}\left(  x^{-1/\gamma_{1}+\left(  1-\eta\right)
/\gamma\pm\epsilon}\right)  .
\]
Consequently, we have%
\[
x^{1/\gamma_{1}-\left(  1-\eta_{0}\right)  /\gamma}\left\{  \mathbf{D}%
_{n}^{\left(  \mathbf{LB}\right)  }\left(  x\right)  -\mathbf{\Gamma}\left(
x;\mathbf{W}\right)  -x^{-1/\gamma_{1}}\dfrac{x^{\tau_{1}/\gamma_{1}}%
-1}{\gamma_{1}\tau_{1}}\sqrt{k}\mathbf{A}_{0}\left(  n/k\right)  \right\}
=o_{\mathbf{P}}\left(  1\right)  ,
\]
uniformly over $x\geq x_{0.}$ Recall that $1/\gamma_{1}=1/\gamma-1/\gamma
_{2},$ then letting $\eta_{0}:=1/2-\xi$ yields $0<\xi<1/2-\gamma/\gamma_{2}$
and achieves the proof.

\subsection{Proof of Theorem \ref{Theorem3}}

The proof is similar, mutatis mutandis, as that of Corollary 3.1 in
\cite{BchMN-16a}. Therefore we omit the details.\medskip

\noindent\textbf{Concluding note}\medskip

\noindent On the basis of Lynden-Bell integration, we introduced a new
estimator for the tail index of right-truncated heavy-tailed data by
considering a random threshold. This estimator may be an alternative to that,
based on Woodroofe integration, recently proposed by \cite{BchMN-16a}. Indeed,
the simulation results show that there is an equivalence between the
asymptotic behaviors of both estimators with respect to biases and rmse's.
However, from a theoritical point of view, the asymptotic normality of the
former requires an additional condition on the sample fraction $k$ of upper
order statistics, namely $\left(  \ref{add}\right)  ,$ which is stonger than
the usual assumption in the context of extremes $\left(  k/n\rightarrow
0\right)  .$

\end{document}